\documentclass[pdflatex,sn-mathphys-num]{sn-jnl}

\usepackage{amsmath,amssymb,amsfonts,amsthm}
\usepackage{mathtools}
\usepackage{booktabs}
\usepackage{array}
\usepackage{float}
\usepackage{microtype}
\usepackage{enumitem}
\usepackage{placeins}
\usepackage[caption=false]{subfig}
\usepackage{braket}
\usepackage[dvipsnames]{xcolor}
\usepackage{tikz}

\hypersetup{
  pdftitle={Quadratic Convexification with Two Switching Variables via Extreme Rays},
  pdfauthor={Yipeng Zhang, Yuyuan Ouyang, Boshi Yang}
}

\renewcommand{\Set}[1]{\left\{#1\right\}}
\setenumerate{label=(\roman*),itemsep=3pt,topsep=3pt}
\mathtoolsset{showonlyrefs}
\allowdisplaybreaks
\newtheorem{theorem}{Theorem}[section]
\newtheorem{proposition}[theorem]{Proposition}
\newtheorem{lemma}[theorem]{Lemma}
\newtheorem{corollary}[theorem]{Corollary}
\theoremstyle{definition}

\theoremstyle{remark}

\DeclareMathOperator{\conv}{conv}
\DeclareMathOperator{\cone}{cone}

\newcommand{\R}{\mathbb R}
\newcommand{\Smat}{\mathbb S}
\newcommand{\vx}{\mathbf{x}}
\newcommand{\vy}{\mathbf{y}}
\newcommand{\vr}{\mathbf{r}}

\newcommand{\vg}{\boldsymbol{\gamma}}

\DeclareRobustCommand{\Sset}{\mathcal{S}}

\newcommand{\Ccone}{\mathcal C}
\newcommand{\Ncone}{\mathcal N}

\begin{document}

\title[Quadratic convexification with two switching variables]
{Quadratic Convexification with Two Switching Variables via Extreme Rays}
\author[1]{\fnm{Yipeng} \sur{Zhang}}
\author[1]{\fnm{Yuyuan} \sur{Ouyang}}
\author[1]{\fnm{Boshi} \sur{Yang}}
\affil[1]{\orgdiv{School of Mathematical and Statistical Sciences}, \orgname{Clemson University}, \orgaddress{\city{Clemson}, \postcode{29634}, \state{South Carolina}, \country{USA}}}

\abstract{%
We consider the convex hull of the quadratic lift \((\vx,\vx\vx^\top,\vy)\) of the set \(\{(\vx,\vy)\mid 0\le\vx\le\vy,\ \vy\in\{0,1\}^n\}\). Descriptions of this convex hull provide convex relaxations for quadratic optimization with switching variables, including problems with indefinite quadratic objectives and constraints. For \(n=2\), we completely characterize the extreme rays of the corresponding cone of nonnegative quadratic functions, which consists of two switching rays, six boundary-product rays, discounted affine-square rays, and a parameterized family of discounted rank-two rays. These families generate the cone, so their linearizations give a complete valid-inequality description of the convex hull in the lifted space. The discounted rank-two rays also show that rank-one liftings cannot generate all valid inequalities. Using this classification, we prove a conjecture of Anstreicher and Burer on an exact semidefinite representation of this convex hull.
}

\maketitle

\section{Introduction}
\label{sec:intro}

Binary variables are frequently used to indicate whether continuous decisions are active. When a binary variable \(y_i\) is coupled with a continuous variable \(x_i\), the condition \(y_i=0\) may be used to force \(x_i=0\), while \(y_i=1\) allows \(x_i\) to vary over a prescribed range. Such binary variables are commonly called \emph{indicator} or \emph{switching variables}. These structures arise, for example, in portfolio optimization \cite{bienstock1996computational}, unit commitment \cite{frangioni2006solving}, and chemical pooling \cite{d2011valid}. For more applications, see \cite{frangioni2006perspective,shafiee2024constrained} and the references therein.

In this paper, we study the quadratic convexification of the fundamental set
\begin{equation}\label{eq:set-n}
\Sset^n:=\Set{(\vx,\vy)\in\R^n\times\{0,1\}^n
\mid 0\leq x_i\leq y_i,\ i=1,\ldots,n}
\end{equation}
and the associated lifted convex hull
\begin{equation}\label{eq:hull-n}
\mathcal K_n:=\conv\Set{(\vx,\vx\vx^\top,\vy)
\mid(\vx,\vy)\in\Sset^n}.
\end{equation}

This structure is closely related to mixed-integer quadratic programs
with indicator variables. Following Dong and Linderoth~\cite{dong2013valid},
consider, for example,
\begin{equation}\label{eq:miqp}
\min\Set{\vx^\top Q\vx+\mathbf c^\top\vx+\mathbf d^\top\vy
\mid (\vx,\vy)\in\Sset^n,\ A\vx+B\vy\leq\mathbf b},
\end{equation}
where \(Q\in\Smat^n\) need not be positive semidefinite.
Here $\Smat^n$ denotes the space of real symmetric $n\times n$ matrices.
We use $A\mathbin{\bullet}B:=\operatorname{tr}(A^\top B)$ for the
Frobenius inner product of real matrices of the same size.
General positive upper bounds can be reduced to the unit bounds in
\(\Sset^n\) by scaling. The objective function can be linearized to 
\(Q\mathbin{\bullet}X+\mathbf c^\top\vx+\mathbf d^\top\vy\) by introducing \(X=\vx\vx^\top\).
Therefore, a description of \(\mathcal K_n\), combined with the additional linear
constraints, yields a convex relaxation of \eqref{eq:miqp}. When there are no additional constraints, the reformulation is exact. The same construction also linearizes quadratic constraints in \(\vx\). In particular, \(\mathcal K_n\) captures the interaction
between the switching variables and all continuous quadratic monomials.

Much of the literature on indicator variables develops convex formulations for epigraphs of specified functions. Frangioni and Gentile~\cite{frangioni2006perspective} derived perspective cuts, and G\"unl\"uk and Linderoth~\cite{gunluk2010perspective} studied perspective formulations for separable convex functions with indicator variables. Subsequent work addresses rank-one and low-rank functions, combinatorial restrictions on the indicators,
and sign restrictions on the continuous variables \cite{wei2022ideal,shafiee2024constrained,han2025compact}. For nonseparable quadratic functions, strong formulations have been obtained by exploiting M-matrix structure
\cite{atamturk2018strong}, positive-semidefinite decompositions \cite{frangioni2020decompositions}, and conic quadratic structure \cite{gomez2021strong}. Related work studies bivariate convex quadratic epigraphs~\cite{han20232} and quadratic constraints with indicator variables~\cite{gomez2024note}. These works generally fix a convex or specially structured function and derive convex hull descriptions or strong formulations for its epigraph or the associated optimization problem.

Another approach is to convexify the quadratic lift \((\vx,\vx\vx^\top,\vy)\), which simultaneously supports every quadratic function of \(\vx\), including indefinite ones.
Dong and Linderoth~\cite{dong2013valid} studied valid inequalities for \(\mathcal K_n\) by lifting valid inequalities for the continuous quadratic hull over \([0,1]^n\), denoted
by \(\mathrm{QPB}_n\) \cite{burer2009nonconvex}. They showed that their lifted-concave-QPB inequalities include the perspective cuts and describe
\[
\mathcal K_n+\Set{(\mathbf 0,P,\mathbf 0) \mid P\in\Smat^n,\ P\succeq0},
\]
the set obtained by allowing an arbitrary positive-semidefinite increment in \(X\).
For convex quadratic minimization, this enlargement preserves the optimal value over \(\mathcal K_n\), since a positive-semidefinite increment in \(X\) cannot decrease the linearized objective. They also observed that concave-QPB inequalities with quadratic coefficient matrices of rank at most one suffice in the continuous setting, and asked to what extent this
remains true after introducing indicators.

For \(n=2\), Anstreicher and Burer~\cite{AnstreicherBurer2010} showed that \(\mathrm{QPB}_2\) is described exactly by the standard positive-semidefinite constraint and the reformulation-linearization technique (RLT) inequalities. For the switching set \(\Sset:=\Sset^2\), they later obtained an exact extended semidefinite formulation for the convex hull that also retains the binary product \(y_1y_2\) \cite{AnstreicherBurerSwitching2021}. For the projected hull \(\mathcal K:=\mathcal K_2\),
denoted by \(\mathcal H'\) in that paper, they conjectured that the standard Shor semidefinite relaxation, together with RLT constraints and the linearized binary and
switching identities, is already exact. The proposed formulation uses one positive-semidefinite constraint of order five and linear constraints.

A closely related result is due to De Rosa and Khajavirad~\cite{de2026explicit}. Whereas the constraints
\(0\leq x_i\leq y_i\) in $\Sset$ impose the upper bounds \(x_i\leq1\), they consider the corresponding unbounded bivariate indicator model
\[
    \{(\vx,\vx\vx^\top,\vy) \mid x_i(1-y_i)=0, \ \vx \geq 0, \ \vy \in \{0,1\}^2\}.
\]
They derive an explicit piecewise description of the closure of its convex hull in the original lifted space and develop a separation procedure based on that description. Their analysis starts from a four-way disjunctive extended formulation. The main difficulty is to project out the auxiliary variables and obtain an explicit description despite the nonlinear constraints. A key insight in their analysis is to recast the projection as a parametric convex optimization problem, which they solve analytically over different regions of the projected space.

We approach these problems from the dual point of view.
Rather than projecting a primal disjunctive formulation, we study the cone of quadratic inequalities valid on \(\Sset\). For \(R\in\Smat^2\), \(\vr,\vg\in\R^2\), and \(p\in\R\), let
\begin{align}
	\label{eq:polynomial}
    q_{R,\vr,p,\vg}(\vx,\vy) := \vx^\top R\vx+2\vr^\top\vx+p+\vg^\top\vy
\end{align}
and define
\begin{equation}\label{eq:cone}
    \Ncone(\Sset) := \{(R,\vr,p,\vg) \mid q_{R,\vr,p,\vg}(\vx,\vy) \geq 0 \quad \forall \ (\vx,\vy) \in \Sset\}.
\end{equation}
We also write $q_{R,\vr,p,\vg}\in \Ncone(\Sset)$, or simply $q\in \Ncone(\Sset)$, to identify a polynomial with its coefficient tuple $(R,\vr,p,\vg)$. Each member of \(\Ncone(\Sset)\) induces the valid lifted inequality
\[
    R \mathbin{\bullet} X + 2\vr^\top\vx + p +\vg^\top\vy \geq 0
\]
for \(\mathcal K\), and every affine inequality valid for \(\mathcal K\) arises in this way. Thus, classifying the extreme rays of \(\Ncone(\Sset)\) gives the supporting inequalities of \(\mathcal K\). The four assignments of \(\vy\in\{0,1\}^2\) reduce nonnegativity on $\Sset$ to simultaneous nonnegativity at the origin, on two unit intervals, and on the unit square. Our analysis exploits the zeros of these four coupled restrictions and their first-order behavior at those points.

Our first contribution is a complete characterization of the extreme rays of \(\Ncone(\Sset)\). They consist of two switching rays, six boundary-product rays, a family of discounted affine-square rays, and a parameterized family of discounted rank-two rays. By analyzing the compatible zero patterns of the four restrictions described above, we prove that these classes exhaust all extreme rays and generate \(\Ncone(\Sset)\). The
classification therefore gives a complete supporting-inequality description of \(\mathcal K\) in the original lifted space \((\vx,X,\vy)\).

Our second contribution establishes the necessity of rank-two lifting
already for $n=2$. The discounted rank-two rays have positive-definite
continuous quadratic coefficient matrices and, by extremality, cannot be
expressed as conic combinations of valid liftings whose continuous quadratic
matrices have rank at most one. This addresses the question raised by
Dong and Linderoth~\cite{dong2013valid} about the sufficiency of rank-one
lifting after introducing indicators.

Our third contribution is a proof of the conjecture of Anstreicher and Burer~\cite{AnstreicherBurerSwitching2021}. We show that the linear and semidefinite constraints in the proposed formulation imply the lifted inequality associated with every classified extreme ray. The complete valid-inequality description then establishes that the formulation is exact.

The remainder of the paper is organized as follows. Section~\ref{sec:switching-classification} characterizes the extreme rays of \(\Ncone(\Sset)\), proves that they generate the cone, and explains the necessity of rank-two lifting. Section~\ref{sec:switching-conjecture} uses this classification to prove the conjectured convex-hull description. Section~\ref{sec:conclusion} summarizes the results and discusses directions for future work.

\section{Extreme-ray classification}
\label{sec:switching-classification}

In this section, we characterize the cone \(\Ncone(\Sset)\) by classifying its extreme rays into four families. The following theorem gives the complete classification.

\begin{theorem}
\label{thm:switching-main}
	The extreme rays of \(\Ncone(\Sset)\) are precisely the rays spanned by the following functions. 
	\begin{enumerate}
		\item The switching rays: \(1-y_1\) and \(1-y_2\).
		\item The boundary-product rays: \(\ell_i\ell_j\),
		\(1\leq i<j\leq4\), where \(\ell_1(\vx):=x_1\),
		\(\ell_2(\vx):=1-x_1\), \(\ell_3(\vx):=x_2\), and
		\(\ell_4(\vx):=1-x_2\).
		\item The discounted affine-square rays:
		\[
		h_\ell(\vx,\vy)
		:=\ell(\vx)^2
		-(1-y_1)\min_{t\in[0,1]}\ell(0,t)^2
		-(1-y_2)\min_{t\in[0,1]}\ell(t,0)^2,
		\]
		where \(\ell\) is a nonzero affine function whose zero set contains
		two distinct points of \([0,1]^2\).
		\item The discounted rank-two rays:
        \begin{equation}\label{eq:rank-two-cut}
        \begin{aligned}
        q_{u,v,a,b}(\vx,\vy)
        ={}&\left(\vx-\binom{u}{v}\right)^\top
        \begin{pmatrix}bv^2&abuv\\abuv&au^2\end{pmatrix}\left(\vx-\binom{u}{v}\right)\\
        &-u^2v^2(1-ab)\bigl[b(1-y_1)+a(1-y_2)\bigr],
        \end{aligned}
        \end{equation}
        where
        the parameters satisfy
        \begin{equation}\label{eq:rank-two-domain}
        a,b>0,\ ab<1,\ 
        0<u\le\frac{1}{1+a},\text{ and } 0<v\le\frac{1}{1+b}.
        \end{equation}
	\end{enumerate}
    Moreover, \(\Ncone(\Sset)\) is the conic hull of its extreme rays.
\end{theorem}

The switching rays depend only on the switching variables, and the
boundary-product rays are products of two distinct affine factors defining
the sides of the unit square. The restrictions of the last two families to
\(\vy=(1,1)\) are, respectively, affine squares and shifted positive-definite
quadratic forms, whose Hessians with respect to \(\vx\) have ranks one and
two. In both families, \emph{discounted} refers to subtracting the minimum
of this restriction on each edge \(x_i=0\), multiplied by \(1-y_i\).

The following subsections develop the ingredients for the proof of
Theorem~\ref{thm:switching-main}, which is completed in
Section~\ref{sec:switching-unique-zero}.
In Section~\ref{sec:switching-preliminaries}, we prove the validity and extremality of the first three ray families. In Section \ref{sec:switching-zero-reduction}, we show that a nonnegative switching polynomial belongs to the cone generated by these families if either interval restriction vanishes identically or its square restriction has at least two distinct zeros. In Section \ref{sec:switching-unique-zero}, we analyze the remaining case and complete the proof.

\subsection{Validity and extremality of the first three ray families}
\label{sec:switching-preliminaries}

The four assignments of the switching variables restrict $\Sset$ to
the origin, two unit intervals,
and the unit square.
For $q=q_{R,\vr,p,\vg}$ as in \eqref{eq:polynomial}, define
\begin{equation}\label{eq:restrictions}
\begin{aligned}
p_q&:=q((0,0),(0,0))=p,\\
f_q(t)&:=q((t,0),(1,0))
       =R_{11}t^2+2r_1t+p+\gamma_1,\\
g_q(t)&:=q((0,t),(0,1))
       =R_{22}t^2+2r_2t+p+\gamma_2,\\
F_q(\vx)&:=q(\vx,(1,1))
       =\vx^\top R\vx+2\vr^\top\vx+p+\gamma_1+\gamma_2,
\end{aligned}
\end{equation}
where $t\in[0,1]$ and $\vx\in[0,1]^2$.
Then $q\in\Ncone(\Sset)$ if and only if all four restrictions
are nonnegative on their respective domains.
These restrictions satisfy
\begin{equation}\label{eq:restriction-relations}
\begin{aligned}
q(\vx,\vy)&=F_q(\vx)-\gamma_1(1-y_1)-\gamma_2(1-y_2),\\
F_q(x_1,x_2)&=f_q(x_1)+g_q(x_2)-p_q+2R_{12}x_1x_2,\\
g_q(t)&=F_q(0,t)-\gamma_1,\\
f_q(t)&=F_q(t,0)-\gamma_2,\\
p_q&=F_q(0,0)-\gamma_1-\gamma_2.
\end{aligned}
\end{equation}
In particular, the four nonnegativity conditions can be written as 
\begin{equation}\label{eq:discount-inequalities}
\begin{aligned}
F_q(\vx)&\ge0 && \forall \, \vx\in[0,1]^2,\\
F_q(0,t)&\ge\gamma_1 &&\forall \, t\in[0,1],\\
F_q(t,0)&\ge\gamma_2 &&\forall \, t\in[0,1],\\
F_q(0,0)&\ge\gamma_1+\gamma_2.
\end{aligned}
\end{equation}

Let $\Ncone([0,1]^2)$ denote the cone of bivariate polynomials of
total degree at most two that are nonnegative on the unit square.
We first characterize its extreme rays, then use them to construct
the boundary-product and discounted affine-square families in
$\Ncone(\Sset)$.

\begin{samepage}
\begin{proposition}\label{prop:switching-square-cone}
The extreme rays of $\Ncone([0,1]^2)$ are precisely those spanned
by the six boundary products $\ell_i\ell_j$, $1\le i<j\le4$,
defined in Theorem~\ref{thm:switching-main}, and by the affine
squares $\ell^2$, where $\ell$ is a nonzero affine function whose
zero set contains two distinct points of $[0,1]^2$.
\end{proposition}
\end{samepage}

\begin{proof}
Let $P\in\Ncone([0,1]^2)$, and write $\mathbf{z}=(1,x_1,x_2)^\top$.
Choose $C,B_{ij}\in\Smat^3$ such that
\[
P(\vx)=\mathbf{z}^\top C\mathbf{z} \quad \text{ and } \quad 
\ell_i(\vx)\ell_j(\vx)=\mathbf{z}^\top B_{ij}\mathbf{z}
\quad(1\le i<j\le4),
\]
and set $E=\operatorname{diag}(1,0,0)$.
The quadratic convex hull of $[0,1]^2$ has an exact description using
positive semidefinite (PSD) constraints and constraints from the
reformulation-linearization technique (RLT)~\cite[Theorem~2]{AnstreicherBurer2010}.
Thus $\mu:=\min_{\vx\in[0,1]^2}P(\vx)$ is the optimum of
\[
\begin{aligned}
\min\quad & C\mathbin{\bullet}Y\\
\text{s.t.}\quad &E\mathbin{\bullet}Y=1,\\
&B_{ij}\mathbin{\bullet}Y\ge0,\qquad 1\le i<j\le4 \\
& Y\succeq0.
\end{aligned}
\]
This semidefinite program is strictly feasible, for example at
\[
Y_0 = \begin{pmatrix}
1&1/2&1/2\\
1/2&1/3&1/4\\
1/2&1/4&1/3
\end{pmatrix}.
\]
Since $\mu$ is finite, strong duality implies that there exist $S\succeq0$ and $\lambda_{ij}\ge0$ such that
\[
C=\mu E+S+\sum_{i<j}\lambda_{ij}B_{ij}.
\]
Therefore,
\[
P(\vx)=\mu+\mathbf{z}^\top S\mathbf{z}
+\sum_{i<j}\lambda_{ij}\ell_i(\vx)\ell_j(\vx).
\]
Since $\mu\ge0$ and $S\succeq0$, the first two terms are a sum
of squares of affine functions.
Thus every extreme ray of $\Ncone([0,1]^2)$ must be spanned by
a boundary product or an affine square.

Suppose that an affine square $\ell^2$ spans an extreme ray.
If $\ell$ has no zero on the square, then $\ell^2$ is strictly
positive there. By compactness, both
$\ell^2+\varepsilon x_1$ and $\ell^2-\varepsilon x_1$
are nonnegative on the square for all sufficiently small
$\varepsilon>0$. Since $x_1$ is not proportional to $\ell^2$,
this contradicts extremality.
If $\ell$ has exactly one zero on the square, that zero is a vertex.
After reflecting coordinates and changing the sign of $\ell$
if necessary, we may write $\ell=ax_1+bx_2$ with $a,b>0$.
Then
\[
\ell^2=a^2x_1^2+2abx_1x_2+b^2x_2^2
\]
is a sum of nonproportional nonnegative quadratics, again
contradicting extremality.
Therefore, the zero set of $\ell$ must contain two distinct
points of the square.

It remains to prove extremality of the stated polynomials,
which are all nonnegative on the square.
Let $P$ be one of these polynomials, and suppose that
$P=P_1+P_2$ with $P_1,P_2\in\Ncone([0,1]^2)$.
Then $0\le P_1\le P$ on $[0,1]^2$.
If $P=\ell_i\ell_j$, then $P_1$ vanishes on the two sides
where $\ell_i$ and $\ell_j$ vanish, respectively.
Hence both affine factors divide $P_1$, so
$P_1=\lambda P$ for some $\lambda\in\mathbb R$.

If $P=\ell^2$, the segment joining two distinct zeros of $\ell$
lies in the square and in the zero set of $\ell$.
Since $P_1$ vanishes on this segment, we may write
$P_1=\ell \hat\ell$ for an affine function $\hat\ell$.
At every point $\vx \in [0,1]^2$ where $\ell(\vx) \ne 0$,
\[
|\hat\ell(\vx)|=\frac{|P_1(\vx)|}{|\ell(\vx)|}\le|\ell(\vx)|.
\]
By continuity, $\hat\ell$ also vanishes on the segment.
Since $\hat\ell$ is affine, $\hat\ell=\lambda\ell$ for some
$\lambda\in\mathbb R$, and again $P_1=\lambda P$.
In either case, $0\le P_1\le P$ implies $0\le\lambda\le1$,
and $P_2=(1-\lambda)P$. Thus $P$ spans an extreme ray.
\end{proof}

We now lift the quadratic functions in $\Ncone([0,1]^2)$ to $\Ncone(\Sset)$ by choosing the
coefficients of the switching variables, following the viewpoint
of Dong and Linderoth~\cite{dong2013valid}.
For a fixed $F\in\Ncone([0,1]^2)$, every polynomial of the form
\eqref{eq:polynomial} with square restriction $F$ can be written as
\[
q(\vx,\vy)=F(\vx)-\gamma_1(1-y_1)-\gamma_2(1-y_2).
\]
Define the two edge minima
\begin{equation}\label{eq:edge-minima}
m_1(F):=\min_{t\in[0,1]}F(0,t) \quad \text{ and } \quad 
m_2(F):=\min_{t\in[0,1]}F(t,0).
\end{equation}
By \eqref{eq:discount-inequalities}, $q \in \Ncone(\Sset)$ if and only if
\begin{equation}\label{eq:valid-discounts}
\gamma_1\le m_1(F),
\gamma_2\le m_2(F),\text{ and } 
\gamma_1+\gamma_2\le F(0,0).
\end{equation}
Increasing either coefficient strengthens the inequality
$q(\vx,\vy)\ge0$ on $0\le\vy\le1$.
The following lemma shows that, when both individual upper bounds
can be attained simultaneously, choosing these coefficients
preserves extremality.

\begin{lemma}\label{lem:switching-greatest-lifting}
Let $F\in\Ncone([0,1]^2)$ satisfy
\begin{equation}\label{eq:discount-compatibility}
m_1(F)+m_2(F)\le F(0,0),
\end{equation}
and define
\[
\mathsf L(F)(\vx,\vy)
:=F(\vx)-m_1(F)(1-y_1)-m_2(F)(1-y_2).
\]
Then $\mathsf L(F)\in\Ncone(\Sset)$, and every
$q\in\Ncone(\Sset)$ with $F_q=F$ admits the decomposition
\begin{equation}\label{eq:lift-decomposition}
q=\mathsf L(F)
+\bigl(m_1(F)-\gamma_1\bigr)(1-y_1)
+\bigl(m_2(F)-\gamma_2\bigr)(1-y_2),
\end{equation}
with nonnegative coefficients.
Moreover, if $F$ spans an extreme ray of $\Ncone([0,1]^2)$,
then a polynomial $q\in\Ncone(\Sset)$ with $F_q=F$ spans an
extreme ray of $\Ncone(\Sset)$ if and only if $q=\mathsf L(F)$.
\end{lemma}

\begin{proof}
By \eqref{eq:valid-discounts} and
\eqref{eq:discount-compatibility},
$\mathsf L(F)\in\Ncone(\Sset)$.
For any $q\in\Ncone(\Sset)$ with $F_q=F$, the same validity
conditions imply $m_i(F)-\gamma_i\ge0$ for $i=1,2$.
The decomposition \eqref{eq:lift-decomposition}
then follows from the first equation in
\eqref{eq:restriction-relations}.

Suppose that $F$ spans an extreme ray of $\Ncone([0,1]^2)$.
To prove extremality of $\mathsf L(F)$, consider a decomposition
$\mathsf L(F)=q^1+q^2$ with $q^1,q^2\in\Ncone(\Sset)$.
Restricting to the square, we have
\[
F=F_{q^1}+F_{q^2}.
\]
By extremality of $F$, there exist $\lambda_1,\lambda_2\ge0$
with $\lambda_1+\lambda_2=1$ such that
$F_{q^k}=\lambda_kF$ for $k=1,2$.
Let $\gamma_i^k$ denote the coefficient of $y_i$ in $q^k$.
By \eqref{eq:valid-discounts}, for each $i,k\in\{1,2\}$,
\[
\gamma_i^k\le m_i(F_{q^k})=\lambda_km_i(F).
\]
For each $i$, we also have
\[
\gamma_i^1+\gamma_i^2
=m_i(F)
=(\lambda_1+\lambda_2)m_i(F).
\]
Therefore, every preceding inequality is an equality:
$\gamma_i^k=\lambda_km_i(F)$ for all $i,k\in\{1,2\}$.
By the first equation in \eqref{eq:restriction-relations},
$q^k=\lambda_k\mathsf L(F)$.
Thus $\mathsf L(F)$ spans an extreme ray of $\Ncone(\Sset)$.

Finally, suppose that $q\in\Ncone(\Sset)$ satisfies $F_q=F$
and $q\ne\mathsf L(F)$.
The decomposition \eqref{eq:lift-decomposition}
expresses $q$ as the sum of $\mathsf L(F)$ and a nonzero
nonnegative combination of $1-y_1$ and $1-y_2$.
The square restrictions of these two summands are $F$ and zero,
respectively. Since $F\ne0$, the summands are nonproportional.
Therefore, $q$ does not span an extreme ray.
\end{proof}

We apply this lemma to the extreme rays identified in
Proposition~\ref{prop:switching-square-cone}.

\begin{corollary}
\label{cor:switching-boundary-square-extremality}
Every polynomial in the boundary-product and discounted
affine-square families of Theorem~\ref{thm:switching-main}
belongs to $\Ncone(\Sset)$ and spans an extreme ray.
\end{corollary}

\begin{proof}
For a boundary product $F=\ell_i\ell_j$, both edges $x_1=0$
and $x_2=0$ contain a zero of $F$.
Thus $m_1(F)=m_2(F)=0$,
\eqref{eq:discount-compatibility} holds,
and $\mathsf L(F)=F$.
The conclusion follows from
Proposition~\ref{prop:switching-square-cone} and
Lemma~\ref{lem:switching-greatest-lifting}.

Now let $F=\ell^2$, where $\ell$ is a nonzero affine function
whose zero set contains two distinct points of $[0,1]^2$.
Since $\mathsf L(F)=h_\ell$, it suffices to verify
\eqref{eq:discount-compatibility}.
Write $\ell=c_0+c_1x_1+c_2x_2$, and assume that $c_0\ge0$ without loss of generality. Set
$d_i:=\min\{0,c_i\}$ for $i=1,2$.
Since $\ell$ has a zero on the square,
\[
0 \geq \min_{\vx\in[0,1]^2}\ell(\vx)=c_0+d_1+d_2.
\]

On the edge $x_1=0$, the affine function $\ell(0,t)=c_0+c_2t$
has minimum $c_0+d_2$.
Since its value at $t=0$ is nonnegative, it has a zero on
$[0,1]$ whenever $c_0+d_2\le0$.
Otherwise, it is strictly positive, and the minimum of its
square is $(c_0+d_2)^2$.
The same argument applies on the edge $x_2=0$.
Thus, with $a_+:=\max\{a,0\}$,
\[
m_1(F)=(c_0+d_2)_+^2,\quad \text{ and } \quad 
m_2(F)=(c_0+d_1)_+^2.
\]

Since $d_1\le0$ and $c_0+d_1+d_2\le0$, we have
\[
(c_0+d_1)_+\le\min\{c_0,-d_2\}
=c_0-(c_0+d_2)_+.
\]
Consequently,
\[
m_1(F)+m_2(F)
\le\bigl((c_0+d_1)_++(c_0+d_2)_+\bigr)^2
\le c_0^2=F(0,0).
\]
Thus \eqref{eq:discount-compatibility} holds,
and $h_\ell$ belongs to $\Ncone(\Sset)$ and spans an extreme ray.
\end{proof}

It remains to verify the switching rays, whose square
restrictions vanish identically.

\begin{lemma}\label{lem:switching-rays-extremality}
The affine functions $1-y_1$ and $1-y_2$ belong to $\Ncone(\Sset)$
and each spans an extreme ray.
\end{lemma}

\begin{proof}
Both polynomials are clearly nonnegative on $\Sset$. We prove that $1-y_1$ spans an extreme ray. The same argument applies to $1-y_2$.

Suppose that
\[
1-y_1=q^1+q^2,
\]
where $q^1,q^2\in\Ncone(\Sset)$.
The square and first interval restrictions of $1-y_1$ vanish
identically, so
\[
F_{q^1}+F_{q^2}\equiv0 \quad \text{ and } \quad 
f_{q^1}+f_{q^2}\equiv0.
\]
Since these restrictions are nonnegative on their respective
domains, we have for each $k=1,2$,
\[
F_{q^k}\equiv0 \quad \text{ and } \quad  f_{q^k}\equiv0.
\]

For $k=1,2$, let $\gamma_i^k$ denote the coefficient of $y_i$
in $q^k$. By \eqref{eq:restriction-relations},
\[
\begin{aligned}
\gamma_1^k&=F_{q^k}(0,0)-g_{q^k}(0)=-g_{q^k}(0)\le0,\\
\gamma_2^k&=F_{q^k}(0,0)-f_{q^k}(0)=0.
\end{aligned}
\]
The first equation in \eqref{eq:restriction-relations} now implies
\[
q^k(\vx,\vy)=F_{q^k}(\vx)-\gamma_1^k(1-y_1)-\gamma_2^k(1-y_2)=g_{q^k}(0)(1-y_1).
\]
Thus both summands are nonnegative multiples of $1-y_1$,
which proves that $1-y_1$ spans an extreme ray.
\end{proof}

These results establish the validity and extremality of the first three
ray families in Theorem~\ref{thm:switching-main}.
The next subsection identifies the zero patterns associated with these
families and prepares the analysis of the last family.

\subsection{Reduction to a unique square zero}
\label{sec:switching-zero-reduction}
Let $\Ccone$ denote the conic hull of the switching, boundary-product,
and discounted affine-square rays in
Theorem~\ref{thm:switching-main}.
Every switching or boundary-product ray except
$(1-x_1)(1-x_2)$ has an interval restriction that vanishes identically.
The square restriction of every boundary-product or discounted
affine-square ray has at least two distinct zeros.
The following lemmas show that either of these zero patterns
forces membership in $\Ccone$.

\begin{lemma}\label{lem:switching-zero-interval}
If \(q\in\Ncone(\Sset)\) and either $f_q\equiv0$ or $g_q\equiv0$, then \(q\in\Ccone\).
\end{lemma}

\begin{proof}
By symmetry, it suffices to study only the case with $f_q\equiv0$. Suppose that $f_q\equiv0$. Then,
\(R_{11}=r_1=0\) and \(\gamma_1=-p\), so
\[
F_q(\mathbf{x})
=
R_{22}x_2^2+2r_2x_2+\gamma_2+2R_{12}x_1x_2.
\]
For each fixed \(x_2\in[0,1]\), the function
\(F_q(\cdot,x_2)\) is affine in $x_1$ with slope \(2R_{12}x_2\). Therefore,
\[
P(x_2)
:=
\min_{x_1\in[0,1]}F_q(x_1,x_2)
=
\begin{cases}
F_q(0,x_2), & R_{12}\ge0,\\
F_q(1,x_2), & R_{12}<0.
\end{cases}
\]
Then $F_q(\mathbf{x})=P(x_2)+B(\mathbf{x})$, where
\[
B(\mathbf{x})
:=
\begin{cases}
2R_{12}x_1x_2, & R_{12}\ge0,\\
-2R_{12}(1-x_1)x_2, & R_{12}<0.
\end{cases}
\]
In particular, \(P\) is a quadratic nonnegative on \([0,1]\),
and \(P(0)=\gamma_2\). Using \eqref{eq:restriction-relations} and \(\gamma_1=-p\), we obtain
\[
q(\mathbf{x},\mathbf{y})
=
p(1-y_1)+B(\mathbf{x})
+P(x_2)-P(0)(1-y_2).
\]
The first two terms belong to \(\mathcal C\), since \(p\ge0\)
and \(B\) is a nonnegative multiple of a boundary product. It remains to show that $P(x_2)-P(0)(1-y_2) \in \mathcal C$.

By the representation of polynomials nonnegative on an interval
\cite[Corollary~II.11.3]{Barvinok2002}, we may write
\[
P(x_2)
=
\lambda x_2(1-x_2)
+\sum_{j=1}^{m}\lambda_j(x_2-a_j)^2,
\]
where $\lambda,\lambda_j\ge0$ and $a_j\in[0,1]$.
Since \(P(0)=\sum_{j=1}^{m}\lambda_j a_j^2\), it follows that
\[
P(x_2)-P(0)(1-y_2) = \lambda x_2(1-x_2) + \sum_{j=1}^{m}\lambda_j \bigl[(x_2-a_j)^2-a_j^2(1-y_2)\bigr].
\]
For each \(j\), the affine function
\(\widetilde{\ell}_j(\mathbf{x})=x_2-a_j\) vanishes on the horizontal
segment \([0,1]\times\{a_j\}\) of the unit square.
Its discounts associated with \(1-y_1\) and \(1-y_2\),
respectively, are
\[
\begin{aligned}
\min_{x_2\in[0,1]}\widetilde{\ell}_j(0,x_2)^2
&=
\min_{x_2\in[0,1]}(x_2-a_j)^2
=0,\\
\min_{x_1\in[0,1]}\widetilde{\ell}_j(x_1,0)^2
&=a_j^2.
\end{aligned}
\]
Thus, each bracketed term is a discounted affine square.
Since \(x_2(1-x_2)\) is a boundary product, the remaining
term also belongs to \(\mathcal C\).
Therefore \(q\in\mathcal C\).
\end{proof}

\begin{lemma}\label{lem:switching-multiple-zeros}
If \(q\in\Ncone(\Sset)\) and $F_q$ has at least two distinct
zeros on \([0,1]^2\), then \(q\in\Ccone\).
\end{lemma}

\begin{proof}
Since \(q\in\Ncone(\Sset)\), the restriction \(F_q\) is nonnegative on \([0,1]^2\). Let $\vx^1, \vx^2 \in [0,1]^2$ be two distinct points such that $F_q(\vx^1) = F_q(\vx^2) = 0$. Choose a nonzero affine function \(\ell\) whose zero line $L:=\{\vx\mid\ell(\vx)=0\}$ passes through $\vx^1$ and $\vx^2$. Let $\mathbf u$ and $\mathbf v$ be the two endpoints of the line segment $L\cap[0,1]^2$, and let
\[
\varphi(t):=F_q((1-t)\mathbf u+t\mathbf v) \quad \text{ for } t\in[0,1]
\]
be the restriction of $F_q$ on $L\cap[0,1]^2$, parameterized by $t$. Then $\varphi$ has degree at most two, is nonnegative on \([0,1]\), and has two distinct zeros. As a result, $\varphi$ is either identically zero or concave with both zeros at endpoints. In either case,
\begin{equation} \label{eq:zero-segment-endpoints}
F_q(\mathbf u)=\varphi(0)=0 \quad \text{ and } \quad 
F_q(\mathbf v)=\varphi(1)=0.  
\end{equation}

Choose the largest multiple of \(\ell^2\) that can be subtracted from \(F_q\) while preserving nonnegativity, i.e., 
\[
\lambda := \max\bigl\{\mu\ge0\mid F_q-\mu\ell^2\ge0\text{ on }[0,1]^2\bigr\}.
\]
The maximum exists because the feasible set is nonempty, closed and bounded. Define 
\[
G:=F_q-\lambda\ell^2.
\]
By construction,
\(G\ge0\) on \([0,1]^2\) and
\(G(\mathbf u)=G(\mathbf v)=0\).
By Proposition~\ref{prop:switching-square-cone},
\[
G=\sum_{k=1}^{N}\widetilde\ell_k^2
  +\sum_{1\le i<j\le4}\beta_{ij}\ell_i\ell_j,
\]
where each $\beta_{ij}\ge0$ and each \(\widetilde\ell_k\) is affine, with its
nonnegative square coefficient absorbed into the factor.
Since every summand is nonnegative on the square,
\[
\widetilde\ell_k(\mathbf u)
=\widetilde\ell_k(\mathbf v)=0,
\]
and hence $\widetilde\ell_k=c_k\ell$ for some $c_k\in\mathbb R$. Consequently,
\[
F_q-\left(\lambda+\sum_{k=1}^{N}c_k^2\right)\ell^2
=\sum_{1\le i<j\le4}\beta_{ij}\ell_i\ell_j
\ge0
\quad\text{on }[0,1]^2.
\]
Maximality of \(\lambda\) therefore gives $c_k = 0$ for each $1 \leq k \leq N$. We conclude that
\[
F_q(\vx)=\lambda\ell(\vx)^2+G(\vx),
\quad \text{where} \ 
G\in\cone\{\ell_i\ell_j\mid 1\le i<j\le4\}.
\]
Thus \(G\) is nonnegative and concave in each coordinate
separately.

Let
\[
m_1:=\min_{x_2\in[0,1]}\ell(0,x_2)^2
\quad \text{ and } \quad
m_2:=\min_{x_1\in[0,1]}\ell(x_1,0)^2.
\]
The corresponding discounted affine-square ray is
\[
h_\ell(\vx,\vy)
=
\ell(\vx)^2-m_1(1-y_1)-m_2(1-y_2).
\]
Using \eqref{eq:restriction-relations}, we obtain
\begin{equation}\label{eq:lift-remainder}
q(\vx,\vy) = \lambda h_\ell(\vx,\vy)+G(\vx)+(\lambda m_1-\gamma_1)(1-y_1)+(\lambda m_2-\gamma_2)(1-y_2).
\end{equation}
To show that the switching coefficients are nonnegative,
we establish
\begin{equation} \label{eq:edge-minima-decomposition}
\min_{x_2\in[0,1]}F_q(0,x_2)=\lambda m_1
\quad \text{ and } \quad
\min_{x_1\in[0,1]}F_q(x_1,0)=\lambda m_2.
\end{equation}
By symmetry, it suffices to prove the first identity.
Since \(G\ge0\) on $[0,1]^2$, 
\[
\min_{x_2\in[0,1]}F_q(0,x_2) = \min_{x_2\in[0,1]} \left\{\lambda \ell(0,x_2)^2 + G(0,x_2) \right\} \ge \lambda m_1.
\]
We show the reverse inequality in three cases.

Case 1: \(L\) meets the left side of the square, i.e., $L \cap (\{0\} \times [0,1]) \neq \emptyset$. In this case, $m_1 = 0$ and at least one endpoint of \(L\cap[0,1]^2\), say $\mathbf u$, lies on that side. Therefore, \eqref{eq:zero-segment-endpoints} implies
\[
\min_{x_2\in[0,1]}F_q(0,x_2) = F_q(\mathbf u) = 0 = \lambda m_1.
\]

Case 2: \(L\) misses the left side and is not the right side, i.e., $L \cap (\{0\} \times [0,1]) = \emptyset$ and $\{1\} \times [0,1] \not\subseteq L$. In this case, \(L\) crosses the interior of the square.
Change the sign of \(\ell\), if necessary, so that \(\ell(0,x_2)>0\) for all \(x_2\in[0,1]\). Since $\ell(0,x_2)$ is affine in $x_2$, its minimum over $[0,1]$ occurs at an endpoint $t_0 \in \{0,1\}$. 
The same \(t_0\) minimizes \(\ell(1,x_2)\), since these
two affine functions have the same slope.
Because \(L\) crosses the interior and \(\ell\) is positive
on the left side, \(\ell\) is negative somewhere on the
right side. Hence
\[
\ell(1,t_0)
=
\min_{x_2\in[0,1]}\ell(1,x_2)
<0<\ell(0,t_0).
\]
Thus \(L\) meets the horizontal side \(x_2=t_0\) at
\((s,t_0)\) for some \(s\in(0,1)\).
This point is an endpoint of \(L\cap[0,1]^2\), so
\(G(s,t_0)=0\). Since $G(x_1, t_0)$ is nonnegative and concave in $x_1$,
\[
0=G(s,t_0) \ge (1-s)G(0,t_0)+sG(1,t_0) \ge 0.
\]
It follows that \(G(0,t_0)=0\). Since $\ell(0,x_2)$ is positive over $[0,1]$, 
\[
\min_{x_2\in[0,1]}F_q(0,x_2) \le F_q(0,t_0)=\lambda\ell(0,t_0)^2=\lambda m_1.
\]

Case 3: \(L\) is the right side, i.e., $\{1\} \times [0,1] \subseteq L$.
In this case \(m_1>0\) and
\(\ell(\vx)^2=m_1(1-x_1)^2\).
Set
\[
\delta:=\min_{x_2\in[0,1]}G(0,x_2)\ge0.
\]
For every \((x_1,x_2)\in[0,1]^2\), concavity in \(x_1\)
and nonnegativity give
\[
G(x_1,x_2) \ge(1-x_1)G(0,x_2)+x_1G(1,x_2) \ge \delta(1-x_1) \ge \delta(1-x_1)^2
=\frac{\delta}{m_1}\ell(x_1,x_2)^2.
\]
Consequently, \(\lambda+\delta/m_1\) is feasible in the
maximization defining \(\lambda\). Maximality implies
\(\delta=0\), and hence
\[
\min_{x_2\in[0,1]}F_q(0,x_2) = \min_{x_2\in[0,1]} \left\{\lambda m_1 + G(0,x_2)\right\}
= \lambda m_1+\delta = \lambda m_1.
\]
This completes the proof of identities \eqref{eq:edge-minima-decomposition}.

Since $g_q(x_2)=F_q(0,x_2)-\gamma_1\ge0$ and $f_q(x_1)=F_q(x_1,0)-\gamma_2\ge0$,
we have \(\lambda m_1 \ge \gamma_1\) and \(\lambda m_2 \ge \gamma_2\).
Every term in \eqref{eq:lift-remainder} therefore
belongs to \(\Ccone\), proving \(q\in\Ccone\).
\end{proof}

Lemmas~\ref{lem:switching-zero-interval} and
\ref{lem:switching-multiple-zeros} show that every
\(q\in\Ncone(\Sset)\setminus\Ccone\) satisfies
\(f_q,g_q\not\equiv0\), and that \(F_q\) has at most one zero
on \([0,1]^2\). The following corollary shows that this zero
exists whenever \(q\) spans an extreme ray.

\begin{corollary}\label{cor:switching-unique-zero}
Let \(q\in\Ncone(\Sset)\setminus\Ccone\) span an extreme ray.
Then \(f_q,g_q\not\equiv0\), and \(F_q\) has exactly one zero
on \([0,1]^2\).
\end{corollary}

\begin{proof}
Set \(c:=\min_{\vx\in[0,1]^2}F_q(\vx)\ge0\). For any \((\vx,\vy)\in\Sset\), if \(\vy \neq (1,1)\), then \(x_1x_2=0\), and
\[
q(\vx,\vy)-cx_1x_2=q(\vx,\vy)\ge0.
\]
If \(\vy=(1,1)\), then \(q(\vx,\vy)=F_q(\vx)\) and \(0\le x_1x_2\le1\), so
\[
q(\vx,\vy)-cx_1x_2 = F_q(\vx)-cx_1x_2 \ge F_q(\vx)-c \ge 0.
\]
Therefore, \(q-cx_1x_2\in\Ncone(\Sset)\), and
\[
q=(q-cx_1x_2)+cx_1x_2
\]
is a decomposition into elements of \(\Ncone(\Sset)\).
By extremality, \(cx_1x_2\) must be a nonnegative multiple of \(q\).
Since \(q\notin\Ccone\) and \(x_1x_2\in\Ccone\), this forces
\(c=0\). Thus \(F_q\) has a zero on \([0,1]^2\).
Lemma~\ref{lem:switching-zero-interval} gives
\(f_q,g_q\not\equiv0\), and Lemma~\ref{lem:switching-multiple-zeros} ensures that the zero
of \(F_q\) is unique.
\end{proof}

\subsection{Stationary zeros and discounted rank-two rays}
\label{sec:switching-unique-zero}

By Corollary~\ref{cor:switching-unique-zero}, every extreme ray outside
\(\Ccone\) has nonzero interval restrictions and a unique zero of its
square restriction. In this subsection, we show that such rays are described exactly by the discounted rank-two ray family in Theorem~\ref{thm:switching-main}. We first characterize two-sided perturbations by
conditions at the zeros, then exclude nonstationary square zeros, and
finally prove that the stationary case yields exactly the discounted rank-two rays.
A zero of an interval restriction is called stationary if its derivative vanishes.
A zero of the square restriction is called stationary if its gradient vanishes.
These definitions also apply at boundary points.

\begin{lemma}\label{lem:switching-perturbation}
Let $q\in\Ncone(\Sset)$ satisfy $f_q,g_q\not\equiv0$, and suppose
that $F_q$ has a unique zero $\mathbf{z}\in[0,1]^2$.
For a polynomial $h$ of the form \eqref{eq:polynomial}, there exists
$\varepsilon>0$ such that
$q\pm\varepsilon h\in\Ncone(\Sset)$ if and only if the following
conditions hold for every $s\in[0,1]$ and $i=1,2$:
\begin{equation}\label{eq:perturbation-conditions}
\begin{aligned}
&p_q=0 \ \Longrightarrow\ p_h=0, && \\
&f_q(s)=0 \ \Longrightarrow\ f_h(s)=0, \quad &&
 f_q(s)=f_q'(s)=0 \ \Longrightarrow\ f_h'(s)=0,\\
&g_q(s)=0 \ \Longrightarrow\ g_h(s)=0, &&
 g_q(s)=g_q'(s)=0 \ \Longrightarrow\ g_h'(s)=0,\\
&F_h(\mathbf{z})=0, &&
 \partial_iF_q(\mathbf{z})=0 \ \Longrightarrow\ \partial_iF_h(\mathbf{z})=0.
\end{aligned}
\end{equation}
Moreover, $q$ spans an extreme ray of $\Ncone(\Sset)$ if and only if
every polynomial $h$ of the form \eqref{eq:polynomial} satisfying
\eqref{eq:perturbation-conditions} is a real multiple of $q$.
\end{lemma}

\begin{proof}
(Necessity)
Suppose that $q\pm\varepsilon h\in\Ncone(\Sset)$ for some
$\varepsilon>0$. Then $\varepsilon|h|\le q$ on $\Sset$, and thus
\[
\varepsilon|p_h|\le p_q,\quad
\varepsilon|f_h|\le f_q,\quad
\varepsilon|g_h|\le g_q,\quad
\varepsilon|F_h|\le F_q
\]
on the corresponding domains. In particular, $p_q=0$ implies
$p_h=0$, and $f_h,g_h,F_h$ vanish at the zeros of
$f_q,g_q,F_q$, respectively.

Both $F_q+\varepsilon F_h$ and $F_q-\varepsilon F_h$ attain
their minimum value $0$ over the unit square at $\mathbf{z}$.
Their partial derivatives satisfy
\[
\partial_i(F_q\pm\varepsilon F_h)(\mathbf{z})
\begin{cases}
\ge0,&z_i=0,\\
=0,&0<z_i<1,\\
\le0,&z_i=1.
\end{cases}
\]
When $\partial_iF_q(\mathbf{z})=0$, these conditions for both signs
force $\partial_iF_h(\mathbf{z})=0$.

At any zero $s$ of $f_q$, both $f_q+\varepsilon f_h$ and
$f_q-\varepsilon f_h$ likewise attain their minimum value $0$
on $[0,1]$. Their ordinary derivatives obey the same
endpoint and interior sign conditions, so $f_h'(s)=0$
whenever $f_q'(s)=0$. The argument for $g_q,g_h$ is identical.
Thus \eqref{eq:perturbation-conditions} holds.

(Sufficiency)
Assume that $h$ satisfies \eqref{eq:perturbation-conditions}.
We first establish a constant $K_F\ge0$ such that $|F_h|\le K_FF_q$ on the unit square.
Let $\mathbf{z}$ be the unique zero of $F_q$, and let $T$ denote the cone of feasible directions at $\mathbf{z}$:
\[
T=\left\{\mathbf{d}\in\mathbb R^2\mid
d_i\ge0\text{ if }z_i=0,\quad
d_i\le0\text{ if }z_i=1,\quad i=1,2
\right\}.
\]
For every $\mathbf{d}\in T$, the point $\mathbf{z}+t\mathbf{d}$ belongs to $[0,1]^2$ for all sufficiently small $t\ge0$.
Since $\mathbf{z}$ minimizes $F_q$ on the square, the first-order necessary optimality condition implies $\nabla F_q(\mathbf{z})^\top \mathbf{d} \geq 0$ for all $\mathbf{d} \in T$. Since each direction $d_i \mathbf{e}_i \in T$, each product $\partial_iF_q(\mathbf{z})d_i$ is nonnegative for $\mathbf{d}\in T$. Here \(\mathbf e_i\) denotes the \(i\)-th standard basis vector. Let
\[
T_0 =\{\mathbf{d}\in T\mid \nabla F_q(\mathbf{z})^\top \mathbf{d}=0\}
=\{\mathbf{d}\in T\mid d_i=0\text{ whenever }\partial_iF_q(\mathbf{z})\ne0\}
\]
be the critical cone at $\mathbf z$ for minimizing $F_q$ over the unit
square, namely the cone of feasible directions along which the first-order
term vanishes.
For every nonzero $\mathbf{d}\in T_0$, uniqueness of the zero and the exact quadratic expansion imply
\[
0<F_q(\mathbf{z}+t\mathbf{d}) = \tfrac12t^2\mathbf{d}^\top\nabla^2F_q\,\mathbf{d}
\]
for all sufficiently small $t>0$. Hence
\[
\mathbf{d}^\top\nabla^2F_q\,\mathbf{d}>0 \qquad \forall \, \mathbf{d}\in T_0\setminus\{0\}.
\]
By \eqref{eq:perturbation-conditions}, $\partial_iF_h(\mathbf{z})=0$ whenever
$\partial_iF_q(\mathbf{z})=0$. Therefore, we may choose $\varepsilon_0>0$ sufficiently small such that both $\nabla(F_q+\varepsilon_0F_h)(\mathbf{z})$ and $\nabla(F_q-\varepsilon_0F_h)(\mathbf{z})$ have the same coordinate signs and zeros as $\nabla F_q(\mathbf{z})$. For each sign,
\[
\nabla(F_q\pm\varepsilon_0F_h)(\mathbf{z})^\top \mathbf{d}\ge0
\qquad \forall \,  \mathbf{d}\in T,
\]
with equality precisely when $\mathbf{d}\in T_0$.

If $T_0\ne\{0\}$, the quadratic form $\mathbf{d}^\top\nabla^2F_q\,\mathbf{d}$ has a positive minimum on the compact set $\{\mathbf{d}\in T_0\mid \|\mathbf{d}\|_2=1\}$, while $|\mathbf{d}^\top\nabla^2F_h\,\mathbf{d}|$ is bounded on the same set. After decreasing $\varepsilon_0$ if necessary, we have
\[
\mathbf{d}^\top\nabla^2(F_q\pm\varepsilon_0F_h)\mathbf{d}>0
\qquad \forall \, \mathbf{d}\in T_0\setminus\{0\}.
\]
When $T_0=\{0\}$, this condition is vacuous.
The second-order sufficient condition on the critical cone $T_0$ for
minimization over a polyhedral set therefore implies that $\mathbf{z}$ is a strict local minimizer of both $F_q+\varepsilon_0F_h$ and $F_q-\varepsilon_0F_h$ over the unit square.

Since both $F_q+\varepsilon_0F_h$ and $F_q-\varepsilon_0F_h$ vanish at $\mathbf{z}$, there exists an open neighborhood $U$ of $\mathbf{z}$, such that both functions are nonnegative on $U \cap [0,1]^2$. Equivalently, $|F_h| \le \varepsilon_0^{-1} F_q$ on $U \cap [0,1]^2$. On the other hand, the compact set $[0,1]^2 \setminus U$ contains no zero of $F_q$, so $|F_h|/F_q$ is continuous and bounded there. Combining the two bounds, we obtain $K_F\ge0$ such that $|F_h|\le K_FF_q$ on $[0,1]^2$.

We now consider the interval restrictions. The nonzero polynomial $f_q$ has finitely many zeros in $[0,1]$. At each such zero $s$, the conditions on $h$ ensure that $f_h(s)=0$.

If $f_q'(s)\ne0$, then $s$ is an endpoint. We may choose $\varepsilon_s>0$ sufficiently small such that $f_q'(s)+\varepsilon_s f_h'(s)$ and $f_q'(s)-\varepsilon_s f_h'(s)$ have the same sign as $f_q'(s)$. Both perturbed functions then increase locally when moving into the interval from $s$, so $s$ is a strict local minimizer of both.
If $f_q'(s)=0$, then $f_h'(s)=0$ as well. Moreover, $f_q''(s)>0$, since $f_q$ is a nonzero quadratic that is nonnegative on $[0,1]$. For sufficiently small $\varepsilon_s>0$,
\[
(f_q\pm\varepsilon_s f_h)'(s)=0 \quad \text{ and } \quad (f_q\pm\varepsilon_s f_h)''(s)>0.
\]
The second-order sufficient condition again establishes a strict local minimum at $s$ for both perturbed functions.

In either case, both local minimum values are zero. There exists an open neighborhood $U_s$ of $s$ such that $|f_h| \le \varepsilon_s^{-1} f_q$ on $U_s \cap [0,1]$. There are finitely many such neighborhoods. On the compact complement of their union, $f_q$ is strictly positive and $|f_h|/f_q$ is bounded. The local bounds and this remaining bound together yield
$K_f\ge0$ such that $|f_h|\le K_f f_q$ on $[0,1]$.
If $f_q$ has no zeros, this conclusion follows directly from compactness. Repeating the argument for $g_q$ and $g_h$, we obtain $K_g\ge0$ with $|g_h|\le K_g g_q$ on $[0,1]$.

Finally, let $K_p$ be $|p_h|/p_q$ if $p_q>0$ and zero otherwise, we have $|p_h|\le K_p p_q$.
Let $\varepsilon = 1/\max\{1,K_F,K_f,K_g,K_p\}$. Then, all four restrictions of $q\pm\varepsilon h$ are nonnegative, so $q\pm\varepsilon h\in\Ncone(\Sset)$.

(Extremality)
Suppose that $q$ spans an extreme ray of $\Ncone(\Sset)$ and $q \pm \varepsilon h \in \Ncone(\Sset)$. Then
\[
q=\tfrac12(q+\varepsilon h)+\tfrac12(q-\varepsilon h)
\]
forces \(h\) to be proportional to \(q\).
Conversely, suppose that every $h$ such that $q \pm \varepsilon h \in \Ncone(\Sset)$ for some $\varepsilon >0$ is a real multiple of $q$. For any decomposition \(q=q^1+q^2\) with \(q^1,q^2\in\Ncone(\Sset)\), set \(h=q^1-q^2\). Then $q+h=2q^1 \in \Ncone(\Sset)$ and $q-h=2q^2 \in \Ncone(\Sset)$. Therefore, $h=\lambda q$ for some $\lambda\in\mathbb R$.
It follows that
 \(q^1=(1+\lambda)q/2\) and \(q^2=(1-\lambda)q/2\).
Since \(q\) is nonzero and both summands are nonnegative on \(\Sset\),
these two coefficients are nonnegative. Hence \(q\) spans an extreme ray.
\end{proof}

We next use the perturbation criterion to exclude nonstationary zeros
of the square restriction.

\begin{lemma}\label{lem:switching-nonstationary-zero}
Let \(q\in\Ncone(\Sset)\setminus\Ccone\) span an extreme ray.
Its square restriction has a unique zero \(\mathbf{z}\), and \(\nabla F_q(\mathbf{z})=0\).
\end{lemma}

\begin{proof}
By Corollary~\ref{cor:switching-unique-zero}, $f_q,g_q\not\equiv0$ and $F_q$ has a unique zero $\mathbf{z}=(z_1,z_2)\in[0,1]^2$. Suppose for a contradiction that $\nabla F_q(\mathbf{z})\ne0$.

Suppose first that one of the partial derivatives is positive. Without loss of generality, assume that $\partial_1F_q(\mathbf{z})>0$. The first-order necessary optimality condition at $\mathbf{z}$ implies $z_1=0$. Take $h=x_1x_2$. Its scalar and interval restrictions vanish identically, and
\[
F_h(\mathbf{z})=\partial_2F_h(\mathbf{z})=0.
\]
Since $\partial_1F_q(\mathbf{z})>0$, all conditions in \eqref{eq:perturbation-conditions} hold. However, $h$ is nonzero and not proportional to $q$, since $f_h\equiv0$ whereas $f_q\not\equiv0$. By Lemma~\ref{lem:switching-perturbation}, this contradicts the extremality of $q$.

As a result, $\nabla F_q(\mathbf{z}) \le 0$, and at least one partial derivative is negative. Without loss of generality, assume that $\partial_1F_q(\mathbf{z})<0$. The first-order necessary optimality condition at $\mathbf{z}$ implies $z_1=1$, so we may write $\mathbf{z}=(1,v)$ with $v\in[0,1]$.

We first show that $g_q$ has a zero on $[0,1]$. Take $h=x_2(1-x_1)$. Its restrictions satisfy
\[
p_h=0,\qquad f_h\equiv0,\qquad g_h(t)=t,\qquad
F_h(\mathbf{z})=\partial_2F_h(\mathbf{z})=0.
\]
If $g_q$ has no zero on $[0,1]$ and $\partial_1F_q(\mathbf{z})<0$, then all conditions in \eqref{eq:perturbation-conditions} hold. The polynomial $h$ is nonzero and not proportional to $q$, because $f_h\equiv0$ whereas $f_q\not\equiv0$. Again, this contradicts the extremality of $q$ by Lemma~\ref{lem:switching-perturbation}. Thus $g_q(s)=0$ for some $s\in[0,1]$.

Next, we show that $f_q$ has a stationary zero on $[0,1]$. Take $h=x_1(1-x_1)$. Its restrictions satisfy
\[
p_h=0,\qquad f_h(t)=t(1-t),\qquad g_h\equiv0,\qquad
F_h(\mathbf{z})=\partial_2F_h(\mathbf{z})=0.
\]
If $f_q$ has no stationary zero on $[0,1]$, then the only zeros $f_q$ may have are at endpoints $0$ and $1$ with multiplicity one. Since $f_h$ vanishes at the endpoints, it is clear that \eqref{eq:perturbation-conditions} holds for this $h$ as well.
Since $h\ne0$ and $g_h\equiv0$ whereas $g_q\not\equiv0$, $h$ is not proportional to $q$, contradicting extremality by Lemma~\ref{lem:switching-perturbation}. Therefore, $f_q(r)=f_q'(r)=0$ for some $r\in[0,1]$.

Since $f_q$ is a nonzero quadratic that is nonnegative on $[0,1]$ and has a stationary zero at $r$, it has the form $f_q(t)=c(t-r)^2$ for some $c>0$. In particular, $f_q'(1)=2c(1-r)\ge0$.
By the second equation in \eqref{eq:restriction-relations},
\[
0>\partial_1F_q(1,v)=f_q'(1)+2R_{12}v,
\]
so $v>0$ and $R_{12}<0$. On the other hand, the same equation, $f_q(r)=g_q(s)=0$, and $p_q\ge0$ imply
\[
0\le F_q(r,s)=2R_{12}rs-p_q\le0.
\]
Therefore, $F_q(r,s)=0$ and $rs=0$. Since $v>0$, the point $(r,s)$ differs from $(1,v)$,
contradicting uniqueness of the zero of $F_q$. Therefore, $\nabla F_q(\mathbf{z})=0$.
\end{proof}

It remains to classify the stationary case. The next proposition proves
both necessity and sufficiency, using the same coefficient comparison
to identify the candidate polynomial and to establish its extremality.

\begin{proposition}\label{prop:switching-stationary-zero}
Let $q$ be a polynomial of the form \eqref{eq:polynomial}. Then $q$ spans an extreme ray of $\Ncone(\Sset)$ and \(F_q\) has a unique zero \(\mathbf{z}\in[0,1]^2\) with \(\nabla F_q(\mathbf{z})=0\), if and only if $q$ is a positive multiple of a polynomial $q_{u,v,a,b}$ in \eqref{eq:rank-two-cut} with parameters satisfying \eqref{eq:rank-two-domain}.
\end{proposition}

\begin{proof}
(Necessity)
Suppose that $q$ spans an extreme ray, and write $\mathbf{z}=(u,v)$.
Since $F_q(u,v)=0$ and $\nabla F_q(u,v)=0$, the exact quadratic
expansion at $(u,v)$ is
\begin{equation}\label{eq:centered}
F_q(\vx)=(\vx-(u,v)^\top)^\top R(\vx-(u,v)^\top).
\end{equation}
In particular, for $t\in[0,1]$,
\[
F_q(t,v)=R_{11}(t-u)^2 \quad \text{ and } \quad
F_q(u,t)=R_{22}(t-v)^2.
\]
Nonnegativity of $F_q$ and uniqueness of its zero imply
$R_{11},R_{22}>0$. Thus $f_q,g_q$ are nonzero strictly convex
quadratics that are nonnegative on $[0,1]$, and each has at most
one zero on $[0,1]$.

We first show that $f_q$ has a zero different from $u$. Take
\[
h(\vx, \vy) =(x_1-u)^2-u^2(1-y_1).
\]
Its restrictions satisfy
\[
p_h=0,\qquad f_h(t)=(t-u)^2,\qquad
g_h\equiv0,\qquad F_h(\vx)=(x_1-u)^2.
\]
If $f_q$ has no zero on $[0,1]$ or $f_q(u)=0$, then all conditions
in \eqref{eq:perturbation-conditions} hold.
Indeed, $F_h(\mathbf{z})=0$ and $\nabla F_h(\mathbf{z})=0$, and
$f_h(u)=f_h'(u)=0$.
However, $h$ is nonzero and not proportional to $q$, since
$g_h\equiv0$ and $g_q\not\equiv0$.
By Lemma~\ref{lem:switching-perturbation}, this contradicts
the extremality of $q$.
Therefore, $f_q$ has a unique zero $r\in[0,1]$ with $r\ne u$.
Applying the same argument to $g_q$, we obtain a unique zero
$s\in[0,1]$ with $s\ne v$.

Since $(r,s)\ne(u,v)$, uniqueness of the square zero and the
second equation in \eqref{eq:restriction-relations} imply
\[
0<F_q(r,s)=2R_{12}rs-p_q.
\]
As $p_q\ge0$, we have $R_{12}>0$ and $r,s>0$.
Moreover, $f_q(u)>0$ and $g_q(v)>0$, since $r\ne u$ and $s\ne v$.
Evaluating the same equation at $(u,v)$, we obtain
\[
p_q=f_q(u)+g_q(v)+2R_{12}uv>0.
\]

We next show that $f_q'(r)=g_q'(s)=0$.
Take
\[
h(\vx, \vy) =(x_1-u)^2-u^2(1-y_1)-(r-u)^2(1-y_2).
\]
Its restrictions satisfy
\[
\begin{aligned}
&p_h=-(r-u)^2, && f_h(t)=(t-u)^2-(r-u)^2,\\
&g_h \equiv 0, && F_h(\vx)=(x_1-u)^2.
\end{aligned}
\]
If $f_q'(r)\ne0$, then the conditions in \eqref{eq:perturbation-conditions} hold because $p_q>0$, $f_h(r)=0$, $g_h\equiv0$, $F_h(\mathbf{z})=0$, and $\nabla F_h(\mathbf{z})=0$.
The polynomial $h$ is again nonzero and not proportional to $q$, since $g_h\equiv0$ and $g_q\not\equiv0$. By Lemma~\ref{lem:switching-perturbation}, this contradicts
the extremality of $q$. Therefore, $f_q'(r)=0$. The same argument applied to $g_q$
shows that $g_q'(s)=0$.

The interval restrictions consequently have the form
\[
f_q(t)=R_{11}(t-r)^2 \quad \text{ and } \quad
g_q(t)=R_{22}(t-s)^2.
\]
Comparing their linear coefficients with those of
\eqref{eq:centered}, using the second equation in
\eqref{eq:restriction-relations}, we obtain
\begin{equation}\label{eq:rank-two-compatibility}
R_{12}v=R_{11}(r-u) \quad \text{ and } \quad
R_{12}u=R_{22}(s-v).
\end{equation}
Since $R_{11},R_{22},R_{12}>0$, $u,v\ge0$, and
$r\ne u$, $s\ne v$, these equations imply
\[
0<u<r\le1 \quad \text{ and } \quad 0<v<s\le1.
\]
In particular, $(u,v)\in(0,1)^2$.
For every nonzero $\mathbf{d}\in\mathbb R^2$, the point $(u,v)+t\mathbf{d}$
belongs to the unit square for all sufficiently small $t>0$.
By uniqueness of the square zero and \eqref{eq:centered},
\[
0<F_q((u,v)+t\mathbf{d})=t^2\mathbf{d}^\top R\mathbf{d}.
\]
Therefore, $R\succ0$.

We now identify $q$ with a member of the discounted rank-two ray family.
Set
\[
a=\frac{r-u}{u} \quad \text{ and } \quad b=\frac{s-v}{v}.
\]
Then $a,b>0$, $r=u(1+a)$, and $s=v(1+b)$.
By \eqref{eq:rank-two-compatibility} and $R\succ0$,
\[
ab=\frac{R_{12}^2}{R_{11}R_{22}}<1.
\]
Together with $r,s\le1$, these relations verify all conditions
in \eqref{eq:rank-two-domain}.
Let $\lambda=R_{12}/(abuv)>0$.
The same compatibility equations imply
\[
R=\lambda
\begin{pmatrix}
bv^2 & abuv\\
abuv & au^2
\end{pmatrix}.
\]
Since $g_q(s)=f_q(r)=0$, \eqref{eq:restriction-relations} implies
\[
q(\vx,\vy)
=F_q(\vx)-F_q(0,s)(1-y_1)-F_q(r,0)(1-y_2).
\]
Using $r=u(1+a)$, $s=v(1+b)$, and the expression for $R$, we obtain
\[
F_q(0,s)=\lambda bu^2v^2(1-ab) \quad \text{ and } \quad
F_q(r,0)=\lambda au^2v^2(1-ab).
\]
Consequently, $q=\lambda q_{u,v,a,b}$, as required.

(Sufficiency)
Take any parameters satisfying \eqref{eq:rank-two-domain},
and set $q=q_{u,v,a,b}$, $r=u(1+a)$, and $s=v(1+b)$.
The matrix $R$ in \eqref{eq:rank-two-cut} has first diagonal
entry $bv^2>0$ and determinant
\[
\det R=abu^2v^2(1-ab)>0,
\]
so $R\succ0$. Hence $F_q$ is nonnegative on the unit square
and has the unique stationary zero $(u,v)$.
The other restrictions are
\begin{equation}\label{eq:rank-two-restrictions}
f_q(t)=bv^2(t-r)^2,\quad
g_q(t)=au^2(t-s)^2,\quad \text{ and } \quad
p_q=abu^2v^2(2+a+b)>0.
\end{equation}
Thus $q\in\Ncone(\Sset)$.
The parameter bounds imply
\[
0<u<r\le1 \quad \text{ and } \quad 0<v<s\le1.
\]
The matrix $R$ and the zeros $r,s$ also satisfy
\eqref{eq:rank-two-compatibility}.

To prove extremality, let $h$ be a polynomial of the form
\eqref{eq:polynomial} satisfying
\eqref{eq:perturbation-conditions}, and let $R_h$
denote its quadratic coefficient matrix.
Since the zeros of $F_q,f_q,g_q$ are stationary, the conditions on
$h$ require
\[
F_h(u,v)=0,\quad \nabla F_h(u,v)=0, \quad f_h(r)=f_h'(r)=0, \quad \text{ and } \quad g_h(s)=g_h'(s)=0.
\]
Therefore,
\[
\begin{aligned}
F_h(\vx)&=(\vx-(u,v)^\top)^\top R_h(\vx-(u,v)^\top),\\
f_h(t)&=(R_h)_{11}(t-r)^2,\qquad
g_h(t)=(R_h)_{22}(t-s)^2.
\end{aligned}
\]
Applying the second equation in \eqref{eq:restriction-relations}
to $h$ and comparing linear coefficients, we obtain
\[
(R_h)_{12}v=(R_h)_{11}(r-u)\quad \text{ and } \quad
(R_h)_{12}u=(R_h)_{22}(s-v).
\]
Since $r-u>0$ and $s-v>0$, comparison with
\eqref{eq:rank-two-compatibility} implies
$R_h=\mu R$, where $\mu=(R_h)_{12}/R_{12}\in\mathbb R$.
It follows that $F_h=\mu F_q$, $f_h=\mu f_q$, and $g_h=\mu g_q$.
Moreover, by \eqref{eq:restriction-relations},
\[
p_h=f_h(0)+g_h(0)-F_h(0,0)=\mu p_q.
\]
Since the four restrictions determine every coefficient of the
polynomial, $h=\mu q$.
By Lemma~\ref{lem:switching-perturbation}, $q$ spans an extreme ray
of $\Ncone(\Sset)$.
\end{proof}

We now combine the validity, extremality, and zero-pattern results
to prove Theorem~\ref{thm:switching-main}.
The final step establishes that the extreme rays generate the entire cone.

\begin{proof}[Proof of Theorem~\ref{thm:switching-main}]
By Lemma~\ref{lem:switching-rays-extremality} and
Corollary~\ref{cor:switching-boundary-square-extremality}, every polynomial
in the first three families belongs to $\Ncone(\Sset)$ and spans
an extreme ray. Proposition~\ref{prop:switching-stationary-zero}
establishes validity and extremality of every polynomial in the
discounted rank-two ray family.

Conversely, let $q\in\Ncone(\Sset)$ span an extreme ray. If $q\in\Ccone$, then $q$ is a positive multiple of a polynomial from one of the first three families.
If $q\notin\Ccone$, by Lemma~\ref{lem:switching-nonstationary-zero}, $F_q$ has exactly
one zero $\mathbf{z}$ on $[0,1]^2$, and $\nabla F_q(\mathbf{z})=0$.
Proposition~\ref{prop:switching-stationary-zero} then implies that
$q$ is a positive multiple of $q_{u,v,a,b}$ for parameters
satisfying \eqref{eq:rank-two-domain}.
Thus the four families account for all extreme rays of
$\Ncone(\Sset)$.

To prove the conic-hull assertion, since $\Ncone(\Sset)$ is a closed convex cone in a finite-dimensional coefficient space, it suffices to show that $\Ncone(\Sset)$ is pointed. Suppose that $q,-q\in\Ncone(\Sset)$. Then $F_q$ vanishes on the unit square and is therefore
the zero polynomial. Moreover, since $f_q(0)=g_q(0)=0$, by \eqref{eq:restriction-relations},
\[
\gamma_1=F_q(0,0)-g_q(0)=0 \quad \text{ and } \quad \gamma_2=F_q(0,0)-f_q(0)=0.
\]
The first equation in \eqref{eq:restriction-relations} now implies
\[
q(\vx,\vy)
=F_q(\vx)-\gamma_1(1-y_1)-\gamma_2(1-y_2)
\equiv0.
\]
Thus $\Ncone(\Sset)$ is pointed, which completes the proof.
\end{proof}

The discounted rank-two rays also clarify the role of rank in indicator
lifting. Dong and Linderoth~\cite{dong2013valid} observed that quadratic
coefficient matrices of rank at most one suffice for continuous
concave-QPB inequalities and asked whether this remains true after
introducing indicators. In our nonnegative inequality convention, the
corresponding liftings have positive-semidefinite continuous quadratic
coefficient matrices. Let $q$ be a discounted rank-two polynomial. Since
$q$ spans an extreme ray of $\Ncone(\Sset)$, any decomposition
$q=\sum_{k=1}^m q^k$ with $q^k\in\Ncone(\Sset)$ forces every nonzero
summand to be a positive multiple of $q$. Each such summand therefore
has a continuous quadratic coefficient matrix of rank two. Thus $q$ cannot be expressed
as a conic combination of valid liftings whose continuous quadratic
coefficient matrices have rank at most one.

In fact, the inequality defined by $q$ is not implied even by the entire (infinite)
family of rank-at-most-one valid liftings. To see this, note that this family is closed and invariant under nonnegative scaling. Since it is contained in the pointed convex
cone $\Ncone(\Sset)$, its conic hull is closed by \cite[Theorem~3.15]{rockafellarwets1998}. The preceding extremality argument therefore places $q$ outside this closed conic hull.
By separation, the linearized inequality defined by $q$ is not implied by this family.

This limitation is reflected in the edge minima. Write
$F_q=L_1^2+L_2^2$ with nonzero affine functions $L_1,L_2$. For $i=1,2$,
\[
m_i(F_q)\ge m_i(L_1^2)+m_i(L_2^2).
\]
Each $L_j$ vanishes at the interior zero of $F_q$, so its discounted
affine-square lifting $h_{L_j}$ is valid. If both inequalities were
equalities, the definitions of the discounts would give
$q=h_{L_1}+h_{L_2}$, contradicting extremality because the summands have
rank-one continuous quadratic coefficient matrices. Therefore, lifting
$F_q$ as a whole yields a strictly larger discount on at least one edge
than lifting its square components separately.

\section{The two-variable switching convex hull}
\label{sec:switching-conjecture}

In this section, we use the extreme-ray classification of $\Ncone(\Sset)$ to prove the conjecture proposed by Anstreicher and Burer in \cite[Section~6]{AnstreicherBurerSwitching2021}. Recall that
\begin{equation}
	\mathcal K :=\conv\Set{(\vx,\vx\vx^\top,\vy)\mid(\vx,\vy)\in\Sset}.
\end{equation}
The conjecture asserts that \(\mathcal K\) is described exactly by the following set:
\begin{align}
\label{eq:conjectured-formulation}
\mathcal K_{\mathrm{SDP}}:=&\left\{(\vx,X,\vy)\ \middle|\
\begin{array}{l}
\vx,\vy\in\R^2,\quad X\in\Smat^2,\\
\exists\ t_1,t_2,T_{12},\alpha_1,\alpha_2,Y_{12}\in\R: \\
t_1=1-y_1, \quad t_2=1-y_2, \\
T_{12}=1+Y_{12}-y_1-y_2, \\
\max\{0,y_1+y_2-1\} \leq Y_{12} \leq \min\{y_1,y_2\}, \\
0 \leq \alpha_1 \leq y_1-Y_{12},\quad 0 \leq \alpha_2 \leq y_2-Y_{12}, \\
X_{12} \geq \max\{0,x_1-\alpha_1+x_2-\alpha_2-Y_{12}\},\\
X_{12} \leq \min\{x_1-\alpha_1,x_2-\alpha_2\},\\
X_{11}\leq x_1\leq y_1,\quad X_{22}\leq x_2\leq y_2, \\
M:= \begin{pmatrix}
1&x_1&x_2&t_1&t_2\\
x_1&X_{11}&X_{12}&0&\alpha_1\\
x_2&X_{12}&X_{22}&\alpha_2&0\\
t_1&0&\alpha_2&t_1&T_{12}\\
t_2&\alpha_1&0&T_{12}&t_2
\end{pmatrix} \succeq 0
\end{array}
\right\}.
\end{align}

As explained by Anstreicher and Burer~\cite[Section~6]{AnstreicherBurerSwitching2021}, the formulation
\eqref{eq:conjectured-formulation} combines the standard
Shor semidefinite relaxation with constraints obtained by RLT.
To see its construction, let $(\vx,\vy)\in\Sset$, and set
\[
\begin{gathered}
X=\vx\vx^\top,\quad t_i=1-y_i\quad(i=1,2),\quad
Y_{12}=y_1y_2,\\
T_{12}=t_1t_2,\quad
\alpha_1=x_1t_2,\quad\alpha_2=x_2t_1.
\end{gathered}
\]
The binary nature of $y_i$ implies $t_i^2=t_i$, while
$0\le x_i\le y_i$ implies $x_it_i=x_i(1-y_i)=0$.
Consequently, the matrix in
\eqref{eq:conjectured-formulation} satisfies
\[
M=\mathbf{z}\mathbf{z}^\top\succeq0 \quad \text{ and } \quad
\mathbf{z}=(1,x_1,x_2,t_1,t_2)^\top.
\]
The Shor relaxation replaces the quadratic products in $\mathbf{z}\mathbf{z}^\top$
by lifted variables and drops the rank-one requirement, retaining
$M\succeq0$ and the linearized consequences of the binary and
complementarity identities. These identities account for the
entries $t_i$ on the diagonal and the zero entries corresponding
to $x_it_i$.

The linear inequalities in \eqref{eq:conjectured-formulation} follow from the original bounds and RLT
constraints. The bounds on $Y_{12}$ are the McCormick inequalities
for the product $y_1y_2$ over $[0,1]^2$.
The bounds on $\alpha_1$ and $\alpha_2$ follow by multiplying
$0\le x_1\le y_1$ by $1-y_2$ and $0\le x_2\le y_2$ by
$1-y_1$, respectively, and linearizing the resulting inequalities:
\[
\begin{aligned}
0\le x_1(1-y_2)\le y_1(1-y_2)
&\quad\longrightarrow\quad 0\le\alpha_1\le y_1-Y_{12},\\
0\le x_2(1-y_1)\le y_2(1-y_1)
&\quad\longrightarrow\quad 0\le\alpha_2\le y_2-Y_{12}.
\end{aligned}
\]
Similarly, the bounds on $X_{12}$ are the linearizations of
\[
x_1x_2\ge0,\quad x_1(y_2-x_2)\ge0,\quad
x_2(y_1-x_1)\ge0,\quad (y_1-x_1)(y_2-x_2)\ge0.
\]
Here the products $x_1y_2$ and $x_2y_1$ are represented by
$x_1-\alpha_1$ and $x_2-\alpha_2$, respectively.
Finally, $X_{ii}\le x_i$ is the linearization of
$x_i(1-x_i)\ge0$, and $x_i\le y_i$ is an original bound.

Thus every $(\vx,\vy)\in\Sset$ admits an extension satisfying
\eqref{eq:conjectured-formulation}, so
$(\vx,\vx\vx^\top,\vy)\in\mathcal K_{\mathrm{SDP}}$.
Since $\mathcal K_{\mathrm{SDP}}$ is convex,
$\mathcal K\subseteq\mathcal K_{\mathrm{SDP}}$.
The conjecture therefore asserts that this standard construction,
with one positive-semidefinite constraint of order five and linear
constraints, already describes the projected hull $\mathcal K$
exactly.

To prove the reverse inclusion, we use the description of
$\Ncone(\Sset)$ from Section~\ref{sec:switching-classification}.
For a polynomial $q=q_{R,\vr,p,\vg}$ of the form \eqref{eq:polynomial}, define
its lifted affine form by
\begin{equation}\label{eq:lifted-form}
\widehat q(\vx,X,\vy)
:=R\mathbin{\bullet}X+2\vr^\top\vx+p+\vg^\top\vy.
\end{equation}
Let $\mathcal E$ denote the collection of polynomials in the four
families of Theorem~\ref{thm:switching-main}.
The next proposition combines separation with the conic-generation
assertion of that theorem.

\begin{proposition}\label{prop:switching-hull-duality}
It holds that
\begin{align*}
\mathcal K
&=\Set{(\vx,X,\vy)\in\R^2\times\Smat^2\times\R^2
\mid\widehat q(\vx,X,\vy)\ge0\quad\forall \, q\in\Ncone(\Sset)}\\
&=\Set{(\vx,X,\vy)\in\R^2\times\Smat^2\times\R^2
\mid\widehat q(\vx,X,\vy)\ge0\quad\forall \, q\in\mathcal E}.
\end{align*}
\end{proposition}

\begin{proof}
Let $\mathcal D$ denote the set on the right-hand side of the first
equality. For every $q\in\Ncone(\Sset)$ and $(\vx,\vy)\in\Sset$,
\[
\widehat q(\vx,\vx\vx^\top,\vy)=q(\vx,\vy)\ge0.
\]
Since $\widehat q$ is affine, this inequality holds on the convex
hull $\mathcal K$. Thus $\mathcal K\subseteq\mathcal D$.

Conversely, suppose that \((\bar\vx, \bar X, \bar\vy) \notin \mathcal K\). Since \(\mathcal K\) is compact and convex, by strong separation,  there exist \(R\in\Smat^2\), \(\vr,\vg\in\R^2\), and \(p\in\R\) such that
\begin{equation} \label{eq:separation}
    R\mathbin{\bullet}X + 2\vr^\top\vx + p + \vg^\top\vy \geq 0 \quad \forall \, (\vx,X,\vy)\in\mathcal K,    
\end{equation}
while $R\mathbin{\bullet}\bar X + 2\vr^\top \bar\vx + p + \vg^\top \bar\vy < 0$. Set $q=q_{R,\vr,p,\vg}$. For every $(\vx,\vy)\in\Sset$, the point $(\vx,\vx\vx^\top,\vy)$ belongs to $\mathcal K$, so \eqref{eq:separation} implies
\[
q(\vx,\vy) =R\mathbin{\bullet}(\vx\vx^\top)+2\vr^\top\vx+p+\vg^\top\vy\ge0.
\]
Hence $q\in\Ncone(\Sset)$, while $\widehat q(\bar\vx,\bar X,\bar\vy)<0$.
Thus $(\bar\vx,\bar X,\bar\vy)\notin\mathcal D$, proving $\mathcal D\subseteq\mathcal K$.

The second equality follows from $\Ncone(\Sset)=\cone(\mathcal E)$ in
Theorem~\ref{thm:switching-main}. Indeed, every lifted inequality
indexed by $\Ncone(\Sset)$ is a nonnegative linear combination
of those indexed by $\mathcal E$, so it suffices to impose
the latter inequalities.
\end{proof}

It remains to verify the extreme-ray inequalities on
$\mathcal K_{\mathrm{SDP}}$. For each parameterized family, we give
an explicit decomposition whose nonnegativity follows from
$M\succeq0$ and the linear constraints.

\begin{proposition}\label{prop:switching-extreme-ray-validity}
Every point of $\mathcal K_{\mathrm{SDP}}$ satisfies
\[
\widehat q(\vx,X,\vy)\ge0\qquad \forall \, q\in\mathcal E.
\]
Consequently, $\mathcal K=\mathcal K_{\mathrm{SDP}}$.
\end{proposition}

\begin{proof}
Fix $(\vx,X,\vy)\in\mathcal K_{\mathrm{SDP}}$ and an extension
satisfying \eqref{eq:conjectured-formulation}.
Set
\[
\beta_i:=y_i-Y_{12}-\alpha_i\qquad(i=1,2).
\]
The defining linear constraints imply
\begin{equation}\label{eq:useful-slacks}
\alpha_i\ge0,\quad\beta_i\ge0\quad(i=1,2),\quad \text{ and } \quad T_{12}\ge0.
\end{equation}
We verify the four families in Theorem~\ref{thm:switching-main}.

We begin with the switching and boundary-product inequalities:
\begin{equation}\label{eq:lifted-elementary-families}
\begin{gathered}
1-y_1\ge0,\quad1-y_2\ge0,\\
x_1-X_{11}\ge0,\quad x_2-X_{22}\ge0,\quad X_{12}\ge0,\\
x_1-X_{12}\ge0,\quad x_2-X_{12}\ge0,\quad
1-x_1-x_2+X_{12}\ge0.
\end{gathered}
\end{equation}
The switching inequalities hold because $1-y_i=t_i$ is a diagonal
entry of $M\succeq0$.
The inequalities in the second line are explicit constraints in
\eqref{eq:conjectured-formulation}.
The bounds $X_{12}\le x_i-\alpha_i$ and $\alpha_i\ge0$ imply
$x_i-X_{12}\ge0$ for $i=1,2$.
Finally,
\[
1-x_1-x_2+X_{12} = \bigl(X_{12}-x_1-x_2+\alpha_1+\alpha_2+Y_{12}\bigr)+\beta_1+\beta_2+T_{12}\ge0.
\]

Next, let $h_\ell$ be a discounted affine-square ray, and write
\[
\ell(\vx)=c_0+c_1x_1+c_2x_2 \quad \text{ and } \quad
\mathbf c=(c_0,c_1,c_2)^\top.
\]
Using the edge minima from \eqref{eq:edge-minima}, set
\[
A:=\sqrt{m_1(\ell^2)} \quad \text{ and } \quad B:=\sqrt{m_2(\ell^2)}.
\]
Its lifted affine form is
\begin{equation}\label{eq:discounted-square-lifted}
\widehat h_\ell(\vx,X,\vy)
=\mathbf c^\top
\begin{pmatrix}1&\vx^\top\\\vx&X\end{pmatrix}
\mathbf c-A^2t_1-B^2t_2.
\end{equation}
Replacing $\ell$ by $-\ell$ does not change $h_\ell$, so we may
assume that $c_0\ge0$.
As in the proof of Corollary~\ref{cor:switching-boundary-square-extremality},
set $d_i=\min\{0,c_i\}$ for $i=1,2$.
Since $\ell$ vanishes on the square, we have $c_0+d_1+d_2\le0$, $A=(c_0+d_2)_+$ and $B=(c_0+d_1)_+$.
Thus $B\le\min\{c_0,-d_2\}=c_0-A$, and hence
\begin{equation}\label{eq:discount-sum}
A+B\le c_0.
\end{equation}

We also have
\begin{equation}\label{eq:edge-discount-inequalities}
\begin{aligned}
A(c_0-A)&\ge0,& A(c_0+c_2-A)&\ge0,\\
B(c_0-B)&\ge0,& B(c_0+c_1-B)&\ge0.
\end{aligned}
\end{equation}
To see the first two inequalities, suppose first that $A>0$.
Then $\ell(0,t)$ has no zero on $[0,1]$ and is positive at $t=0$,
so it is positive throughout the interval.
Consequently,
\[
A=\min_{t\in[0,1]}\ell(0,t)\le\min\{c_0,c_0+c_2\}.
\]
If $A=0$, both inequalities hold trivially.
The last two inequalities follow by the same argument on the edge
$x_2=0$.

To use the constraint $M\succeq0$, consider the affine expression
$\ell(\vx)-At_1-Bt_2$ and its coefficient vector
\[
\boldsymbol\omega=(c_0,c_1,c_2,-A,-B)^\top.
\]
The quadratic form $\boldsymbol\omega^\top M\boldsymbol\omega$
is the linearization of its square, using the identities
$t_i^2=t_i$ and $x_it_i=0$ built into $M$.
Expanding this quadratic form and comparing it with
\eqref{eq:discounted-square-lifted}, we obtain
\[
\begin{aligned}
\widehat h_\ell(\vx,X,\vy)
={}&\boldsymbol\omega^\top M\boldsymbol\omega + 2A\bigl((c_0-A)t_1+c_2\alpha_2\bigr) + 2B\bigl((c_0-B)t_2+c_1\alpha_1\bigr) -2ABT_{12}.
\end{aligned}
\]
Substituting
\[
t_1=\alpha_2+\beta_2+T_{12} \quad \text{ and } \quad
t_2=\alpha_1+\beta_1+T_{12},
\]
and collecting the terms involving $T_{12}$ yields the following
decomposition:
\begin{equation}\label{eq:square-certificate}
\begin{aligned}
\widehat h_\ell(\vx,X,\vy)
={}&\boldsymbol\omega^\top M\boldsymbol\omega+2A(c_0-A)\beta_2+2A(c_0+c_2-A)\alpha_2\\
&+2B(c_0-B)\beta_1+2B(c_0+c_1-B)\alpha_1+2\bigl((A+B)(c_0-A-B)+AB\bigr)T_{12}.
\end{aligned}
\end{equation}
Every term on the right-hand side is nonnegative:
the first by $M\succeq0$, the next four by
\eqref{eq:useful-slacks} and
\eqref{eq:edge-discount-inequalities}, and the last by
\eqref{eq:useful-slacks} and
\eqref{eq:discount-sum}.
Thus $\widehat h_\ell(\vx,X,\vy)\ge0$.

Finally, let $q=q_{u,v,a,b}$ be a discounted rank-two ray with parameters satisfying
\eqref{eq:rank-two-domain}, and set
\[
R=
\begin{pmatrix}
bv^2&abuv\\
abuv&au^2
\end{pmatrix},
\quad \text{ and } \quad
D=
\begin{pmatrix}
-u&1&0&u&-au\\
-v&0&1&-bv&v
\end{pmatrix}.
\]
Expanding $R\mathbin{\bullet}(DMD^\top)$ and using
$t_i=1-y_i$ and the definition of $q_{u,v,a,b}$ in
\eqref{eq:rank-two-cut}, we obtain
\begin{equation}\label{eq:rank-two-certificate}
\widehat q(\vx,X,\vy)
=R\mathbin{\bullet}(DMD^\top)
+2abu^2v^2(1-ab)T_{12}.
\end{equation}
Since $bv^2>0$ and
\[
\det R=abu^2v^2(1-ab)>0,
\]
we have $R\succ0$. Moreover, $M\succeq0$ implies
$DMD^\top\succeq0$, so
$R\mathbin{\bullet}(DMD^\top)\ge0$.
The second term in \eqref{eq:rank-two-certificate} is
nonnegative because $a,b,u,v>0$, $ab<1$, and $T_{12}\ge0$.
Thus $\widehat q(\vx,X,\vy)\ge0$.

We have verified all four families.
Proposition~\ref{prop:switching-hull-duality} now implies
$\mathcal K_{\mathrm{SDP}}\subseteq\mathcal K$.
Together with the reverse inclusion established above, this proves
$\mathcal K=\mathcal K_{\mathrm{SDP}}$.
\end{proof}

\section{Conclusion}
\label{sec:conclusion}

We have characterized the extreme rays of \(\Ncone(\Sset)\)
for the bounded two-variable switching set.
The switching, boundary-product, discounted affine-square,
and discounted rank-two families generate this cone and give a complete
valid-inequality description of \(\mathcal K\).
The extremality of the discounted rank-two rays shows that rank-one
liftings cannot generate all valid inequalities.
Using this classification, we have proved the conjecture of
Anstreicher and Burer: the standard Shor semidefinite relaxation,
together with the specified RLT constraints, describes
\(\mathcal K\) exactly after projection.
The proof connects the explicit inequality families with the
extended semidefinite formulation through nonnegativity
certificates.

Extending our approach to $\mathcal K_3$ presents an additional difficulty.
Proposition~\ref{prop:switching-square-cone} relies on the exact PSD--RLT description of $\mathrm{QPB}_2$ to characterize the extreme rays of $\Ncone([0,1]^2)$.
For $\mathrm{QPB}_3$, PSD--RLT is no longer exact, and the
known disjunctive formulation does not directly yield an
analogous extreme-ray characterization
\cite{AnstreicherBurer2010}.
For larger quadratically constrained quadratic programs
with bounded nonnegative continuous variables and indicators,
the discounted rank-two inequalities can be applied, after scaling,
to selected pairs of variables to strengthen relaxations
based on individual perspective constraints.
Developing efficient separation procedures for this family
is a natural direction for future work.

\section*{Acknowledgment}
The first and second authors are partially supported by AFOSR grant FA9550-25-1-0278. This work was supported in part by OpenAI API credits provided by Clemson University and administered by Clemson University Research Computing and Data (RCD). 

The authors used ChatGPT (OpenAI; GPT-5.6 Sol and subsequent models, last accessed November 2026) throughout the development of this work, including the formulation of results, proof development, and manuscript preparation. The authors independently verified and revised the mathematical arguments and take full responsibility for the content.

\bibliography{article}

\end{document}